\documentclass[11pt,a4paper]{amsart}

\usepackage{amscd,amssymb,amsopn,amsmath,amsthm,graphics,amsfonts,mathrsfs,accents,enumerate,verbatim,calc}
\usepackage{enumitem}
\usepackage[dvips]{graphicx}
\usepackage[colorlinks=true,linkcolor=red,citecolor=blue]{hyperref}
\usepackage{latexsym}
\usepackage{mathrsfs}
\usepackage{newclude}
\usepackage{syntonly}
\usepackage{tikz}		%箭图
\usepackage{tikz-cd}	%箭图
\usepackage[all]{xy}
\usepackage{nicematrix}
\usepackage{arydshln} 	%分块矩阵
\usepackage[normalem]{ulem} %删除线
\usepackage{ctex}

\date{}
\hypersetup
{colorlinks=true,linkcolor=black}
\allowdisplaybreaks[4] \footskip=15pt
\renewcommand{\uppercasenonmath}[1]{}

\numberwithin{equation}{section} \theoremstyle{plain}
\newtheorem{lemm}{Lemma}[section]

\newtheorem{prop}[lemm]{Proposition}
\newtheorem{theo}[lemm]{Theorem}

\newtheorem{exam}[lemm]{Example}

\definecolor{morpink}{RGB}{239,229,248}
\definecolor{morblue}{RGB}{146,172,209}
\definecolor{defcolor}{RGB}{50,205,50}
\definecolor{tcolor}{RGB}{210,105,30} % default: 0,124,53
\definecolor{lcolor}{RGB}{210,105,30} % default: 153,255,153
\definecolor{pcolor}{RGB}{253,249,238} % default: 216,255,216
\definecolor{neww}{RGB}{181, 168, 156} % {0,0,0}

\def\no{\noindent}
\def\proof{\no\rmm{Proof:}}
\def\endd{$\hfill\Box$}  				%结束符号

\newcommand{\rmm}[1]{
	{\rm #1}
}

\newcommand{\Rmnum}[1]{\uppercase\expandafter{\romannumeral#1}}

\newlist{rlist}{enumerate}{1}
\setlist[rlist]{label=\roman*\rmm{)}}

\begin{document}
\begin{center}
{\Large  \bf  Construction of ideal cotorsion pairs via recollements of triangulated categories}

\vspace{0.5cm} Qikai Wang and Haiyan Zhu\footnote{Corresponding author

Supported by the National Natural Science Foundation of China  (12271481).}  \\
School of Mathematical Science, Zhejiang University of Technology, Hangzhou 310023, China

Emails: qkwang@zjut.edu.cn and hyzhu@zjut.edu.cn
\end{center}

\bigskip
\centerline { \bf  Abstract}
\leftskip10true mm \rightskip10true mm \noindent
Let $(\mathcal{T}',\mathcal{T},\mathcal{T}'')$ be a recollement of triangulated categories.
A complete ideal cotorsion pair in $\mathcal{T}$ induces complete ideal cotorsion pairs in $\mathcal{T}'$ and $\mathcal{T}''$.
In addition, if $(\mathcal{I}, \mathcal{I}^\perp )$ and $(\mathcal{J},\mathcal{J}^\perp)$ are two complete ideal cotorsion pairs in a triangulated category, then $(\mathcal{I}\cap\mathcal{J}, \langle\mathcal{I}^\perp,\mathcal{J}^\perp\rangle)$ is also a complete ideal cotorsion pair.
By this method, starting from two complete ideal cotorsion pairs in $\mathcal{T}'$ and $\mathcal{T}''$, one can induce a family of complete ideal cotorsion pairs in $\mathcal{T}$.
\\[2mm]

\leftskip0true mm \rightskip0true mm

%{\footnotesize\tableofcontents\label{contents}}

% !TeX root = ideal.tex
\section{Introduction}

The central role in approximation theory for the case of modules, or more general abelian or exact categories, is played by the notion of cotorsion pair, cf. \cite{RJ}.
Fu, Herzog et al. \cite{Fu2013,Fu2016} extended the study of classical cotorsion pairs to ideal cotorsion pairs, in  exact categories.
Their construction relies on ME‑extensions, which are specific to exact categories.
Furthermore, Breaz and Modoi \cite{icpt} extended ideal approximation theory to triangulated categories, who extended the ME-extensions to construct ideal cotorsion pairs in the triangulated category using the Toda bracket.

The construction of (ideal) cotorsion pairs is a recurring theme in relative homological algebra. 
In exact categories, intersection methods have been used to produce classical and ideal cotorsion pairs; see, for instance, \cite{iccp,fp}. 
Similar ideas also appear in the triangulated setting, e.g.\ in \cite{mtcr}. 
In the present paper we develop an intersection technique for ideal cotorsion pairs in triangulated categories.

\begin{theo}\label{intersection}
	Let $\mathcal{I}_1$ and $\mathcal{I}_2$ be special precovering ideals in $\mathcal{A}$ with respect to $\mathfrak{C}$. Then $\mathcal{I}_1\cap\mathcal{I}_2$ is a special precovering ideal in $\mathcal{A}$.
\end{theo}

Recollements of triangulated categories, introduced by Beilinson, Bernstein and Deligne in \cite{fp}, provide a powerful mechanism for transferring structures between categories. 
In the classical case, Chen \cite{cprt} showed that cotorsion pairs can be constructed along a recollement. 
It is therefore natural to ask whether recollements can also be used to lift (or glue) ideal cotorsion pairs. 
The purpose of this paper is to give an affirmative answer and to describe explicit constructions.

\begin{theo}\label{main1}
	Consider the following recollement of triangulated categories:
	\begin{center}
		\begin{tikzcd}
			\mathcal{T}' \arrow[rr, "i_*"] &  & \mathcal{T} \arrow[rr, "j^*"] \arrow[ll, "i^*"', bend right=49] \arrow[ll, "i_!"', bend left=49] &  & \mathcal{T}'' \arrow[ll, "j_*"', bend left=49] \arrow[ll, "j_!"', bend right=49]
		\end{tikzcd}
	\end{center}
	$(\mathcal{I}',\mathcal{J}')$ is a complete ideal cotorsion pair in $\mathcal{T}'$ and $(\mathcal{I}'',\mathcal{J}'')$ is a complete ideal cotorsion pair in $\mathcal{T}''$ if and only if
	there exists a complete ideal cotorsion pair $(\langle  j_!\mathcal{I}'',i_*\mathcal{I}'\rangle ,\langle  i_*\mathcal{J}',j_*\mathcal{J}''\rangle )$ in $\mathcal{T}$.
\end{theo}

The notations $\langle j_!\mathcal{I}'',i_*\mathcal{I}'\rangle$ and $\langle i_*\mathcal{J}',j_*\mathcal{J}''\rangle$ denote the Toda bracket (see Definition \ref{Toda0}) of $j_!\mathcal{I}''$, $i_*\mathcal{I}'$ and of $i_*\mathcal{J}'$, $j_*\mathcal{J}''$.

In addition, starting from $(\mathcal{I}',\mathcal{J}')$ and $(\mathcal{I}'',\mathcal{J}'')$ we obtain further complete ideal cotorsion pairs in $\mathcal{T}$. 
For the reader's convenience, the relationships among these constructions are summarized in a diagram in the appendix.

The paper is organized as follows. 
In Section~2 we collect the necessary background on ideal cotorsion pairs in triangulated categories, including almost exact structures and phantom ideals, (special) precovering/preenveloping ideals, and Toda brackets, together with several auxiliary results. 
In Section~3 we develop an intersection method for ideals and prove Theorem~\ref{intersection}. 
In Section~4 we apply these preparations to recollements of triangulated categories: we first study the transport of ideal cotorsion pairs along the recollement functors and then establish the gluing result of Theorem~\ref{main1}.

\include*{Pre}
% !TeX root = ideal.tex
\section{Construct Ideal cotorsion pair by intersection}

In this section, we show that the intersection of two special precovering ideals is also a special precovering ideal.
We start with the following lemmas.

\begin{lemm}
	If $\mathcal{I}_1$ and $\mathcal{I}_2$ are ideals in $\mathcal{A}$, then $\mathcal{I}_1\cap\mathcal{I}_2$ is also an ideal in $\mathcal{A}$.
\end{lemm}

\proof
It follows immediately from the definition that $0\in \mathcal{I}_1\cap\mathcal{I}_2$.
Since an ideal is closed under compositions with arbitrary morphisms and finite direct sums, the intersection of ideals is also closed under compositions with arbitrary morphisms and finite direct sums.
\endd

Next, we consider the intersection of object ideals. This is also an important reason for studying the intersection of ideals.

\begin{lemm}
Let $\mathcal{A}$ be an extension-closed subcategory of $\mathcal{T}$ and $\mathfrak{C}$ an almost exact structure on $\mathcal{A}$. If $\mathcal{D}$ and $\mathcal{E}$ are two classes of objects in $\mathcal{A}$ closed under finite direct sums, then $Ideal(\mathcal{E}\cap\mathcal{D})\subseteq Ideal(\mathcal{E}) \cap Ideal(\mathcal{D})$.
\end{lemm}
\proof
For any morphism $i\in Ideal(\mathcal{E}\cap\mathcal{D})$, there exists an object $X\in \mathcal{E}\cap\mathcal{D}$ such that $1_X$ factors through $i$.
Then $i\in Ideal(\mathcal{E})$ and $i\in Ideal(\mathcal{D})$, that is, $i\in Ideal(\mathcal{E}) \cap Ideal(\mathcal{D})$.\endd

\begin{exam}
	Let $\mathcal{D}$ be the class generated by an object $D$, and $\mathcal{E}$ be the class generated by an object $E$ with $D\neq E$.
	Suppose that there exists a nonzero morphism in $\mathrm{Hom}(D,E)$ or $\mathrm{Hom}(E,D)$.
	Then the ideal $Ideal(\mathcal{E}\cap\mathcal{D})$ is empty, but $Ideal(\mathcal{E}) \cap Ideal(\mathcal{D})$ is not.
\end{exam}

The above lemma shows that an intersection of object ideals can contain more elements than the ideal of the intersection of objects, which enables an intersection of ideals to have properties that the ideal of the intersection of objects does not possess.

Next, we show that the intersection of two special precovering ideals is also a special precovering ideal.

\noindent \textit{\textbf{Proof of Theorem \ref{intersection}}}:

Let $A$ be any object in $\mathcal{A}$. Since $\mathcal{I}_1$ and $\mathcal{I}_2$ are special precovering ideals, there exist distinguished triangles for $k=1,2$:
$$ K_k \longrightarrow B_k \xrightarrow{i_k} A \xrightarrow{g_k} K_k[1] $$
such that $i_k \in \mathcal{I}_k$, and the morphisms $g_{k}$ decompose as $g_k = j_k[1] \circ \phi_k$, where $j_k \in \mathcal{I}_k^\perp$ and $\phi_k \in \Phi(\mathfrak{C})$.

Consider the morphisms $\phi_1: A \to D_1[1]$ and $\phi_2: A \to D_2[1]$ in $\Phi(\mathfrak{C})$.
Since $\Phi(\mathfrak{C})$ is closed under finite coproducts, the direct sum morphism $\phi_1 \oplus \phi_2: A \oplus A \to D_1[1] \oplus D_2[1]$ belongs to $\Phi(\mathfrak{C})$.
Since $\Phi(\mathfrak{C})$ is closed under base change, we can pull back along the diagonal morphism $\Delta = (1_A, 1_A)^T: A \to A \oplus A$.
Define $\Psi: A \to D_1[1] \oplus D_2[1]$ as:
$$ \Psi := (\phi_1 \oplus \phi_2) \circ \Delta = \begin{pmatrix} \phi_1 \\ \phi_2 \end{pmatrix}. $$
Thus, $\Psi \in \Phi(\mathfrak{C})$.
Define $J: D_1 \oplus D_2 \to K_1 \oplus K_2$ as the diagonal matrix:
$$ J := \begin{pmatrix} j_1 & 0 \\ 0 & j_2 \end{pmatrix}. $$
Since $j_1 \in \mathcal{I}_1^\perp$ and $j_2 \in \mathcal{I}_2^\perp$, by the properties of orthogonal ideals (specifically that $(\mathcal{I}_1 \cap \mathcal{I}_2)^\perp$ contains sums of morphisms from $\mathcal{I}_1^\perp$ and $\mathcal{I}_2^\perp$), we have $J \in (\mathcal{I}_1 \cap \mathcal{I}_2)^\perp$.

Define the morphism $G: A \to K_1[1] \oplus K_2[1]$ by:
$$ G := J[1] \circ \Psi = \begin{pmatrix} j_1[1] \circ \phi_1 \\ j_2[1] \circ \phi_2 \end{pmatrix} = \begin{pmatrix} g_1 \\ g_2 \end{pmatrix}. $$

Complete $G$ to a distinguished triangle in $\mathfrak{C}$:
$$K_{1}\oplus K_{2}\rightarrow  P \xrightarrow{f} A \xrightarrow{G} K_1[1] \oplus K_2[1]. $$
It remains to show $f \in \mathcal{I}$.
Using the projections $\pi_k: K_1[1] \oplus K_2[1] \to K_k[1]$ for $k=1,2$, we have:
$$ g_k \circ f = \pi_k \circ G \circ f = 0. $$
Since $g_k \circ f = 0$, $f$ factors through $i_k \in \mathcal{I}_k$. Since $\mathcal{I}_k$ is an ideal, this implies $f \in \mathcal{I}_k$.
Therefore, $f \in \mathcal{I}_1 \cap \mathcal{I}_2$.
\endd

% !TeX root = ideal.tex
\section{Ideal cotorsion pair in recollement of triangulated categories}

Note that if $\mathcal{A}=\mathcal{T}$ and $\mathfrak{C}$ is the class of all triangles then every precovering ideal is special. Dually, every preenveloping ideal in $\mathcal{T}$ is special with respect to the class of all triangles in $\mathcal{T}$.
Thus, next we only consider the case where $\mathcal{A}=\mathcal{T}$ and $\mathfrak{C}$ is the class of all triangles.
In this case, the Toda bracket $\langle \mathcal{I},\mathcal{J}\rangle_{\mathfrak{C}}$ is denoted by $\langle \mathcal{I},\mathcal{J}\rangle$ for simplicity.
And the condition of having enough $\mathfrak{C}$-injective morphisms and having enough $\mathfrak{C}$-projective morphisms are satisfied in our setting, since the $\mathfrak{C}$ is the class of all triangles in $\mathcal{A}$.

\subsection{Functorial transport along a recollement}\

In this subsection, we show that an ideal cotorsion pair in $\mathcal{T}$ can induce ideal cotorsion pairs in $\mathcal{T}'$ and $\mathcal{T}''$ under some conditions.

\begin{lemm}\label{main2}
	Consider a recollement as in Definition \ref{recollement}.
	Let $(\mathcal{I},\mathcal{J})$ be a (complete) ideal cotorsion pair in $\mathcal{T}$.
	Then:

	$\mathrm{(1)}$ If $j_!j^*\mathcal{I} \subseteq \mathcal{I}$, then $(j^*\mathcal{I}, j^*\mathcal{J})$ is a (complete) ideal cotorsion pair in $\mathcal{T}''$.

	$\mathrm{(2)}$ If $i_*i^*(\mathcal{I})\subseteq\mathcal{I}$, then $(i^*\mathcal{I}, i^!\mathcal{J})$ is a (complete) ideal cotorsion pair in $\mathcal{T}'$.

\end{lemm}

\proof
We prove the case for $j^*$; the another case follows by duality.

\textbf{Step 1: Stability.}
We deduce the stability of $\mathcal{J}$ from the stability of $\mathcal{I}$.
We show that $j_*j^*(\mathcal{J})\subseteq\mathcal{J}$.
Let $\beta:A\rightarrow B$ be a morphism in $\mathcal{J}$.
To show $j_*j^*(\beta) \in \mathcal{J}$, we consider the composition $\xi = j_*j^*(\beta)[1] \circ \phi \circ \alpha$ with any $\alpha \in \mathcal{I}$ and any phantom $\phi\in \Phi(\mathfrak{C})$.
Consider the adjoint pair $(j^*, j_{*})$.
The morphism $\xi$ corresponds to the morphism $\xi' = j^*(\beta)[1]\circ \epsilon_{j^*A} \circ j^*(\phi)\circ j^*(\alpha)$,
where $\epsilon_{j^*A}: j^*j_{*}(j^*A) \rightarrow j^*A$ is the counit of the adjoint pair $(j^*, j_{*})$.
Furthermore, consdier the adjoint pair $(j_!, j^*)$.
The morphism $\xi'$ corresponds to the morphism $\tilde{\xi} = \beta[1]\circ \epsilon_{A} \circ j_!(\epsilon_{j^*A}) \circ j_!j^*(\phi)\circ j_!j^*(\alpha)$.
Since $j_{!}j^*(\mathcal{I})\subseteq \mathcal{I}$, 
and $(\mathcal{I}, \mathcal{J})$ is an ideal cotorsion pair, such a composition must vanish, so $\tilde{\xi} = 0$.
By the bijection of the adjunction, $\xi = 0$.
Therefore, $j_*j^*(\beta) \in \mathcal{J}$.

\textbf{Step 2: Orthogonality.}
We show that $j^*\mathcal{I} \perp j^*\mathcal{J}$ in $\mathcal{A}''$.
Let $\alpha''$ be a morphism in $j^*\mathcal{I}$, $\beta''$ be a morphism in $j^*\mathcal{J}$ and $\varphi''$ be a phantom morphism in $\mathcal{A}''$.
By the definition of the image ideal, $\alpha''$ can be written as a composition:
$\alpha'' = f'' \circ j^*(\alpha) \circ g''$
where $\alpha: A \to B$ is a morphism in $\mathcal{I}$, and $f'', g''$ are morphisms in $\mathcal{A}''$.

To show $\beta'' \circ \varphi'' \circ \alpha'' = 0$, it suffices to show $\beta'' \circ \varphi'' \circ f'' \circ j^*(\alpha) = 0$.
Consider the adjoint pair $(j^*, j_*)$.
The composition $\Psi'' = \beta'' \circ \varphi'' \circ f'' \circ j^*(\alpha)$ in $\mathcal{A}''$ corresponds to the following morphism in $\mathcal{A}$:
\[ \Psi = j_*(\beta'') \circ j_*(\varphi'') \circ j_*(f'') \circ \eta_{B} \circ \alpha \]
where $\eta_{B}: B \to j_*j^*(B)$ is the unit morphism.

From Step 1, we established that $j_*(\beta'') \in \mathcal{J}$.

Since $(\mathcal{I}, \mathcal{J})$ is a cotorsion pair in $\mathcal{A}$, we have $\Psi = 0$.
By the isomorphism of the adjunction, $\Psi = 0$ implies $\Psi'' = 0$.
Therefore, $\beta'' \circ \varphi'' \circ \alpha'' = \Psi'' \circ g'' = 0$.
Thus, $j^*\mathcal{I} \perp j^*\mathcal{J}$.

\textbf{Step 3: Completeness.}

\textbf{(i)} We show that $(j^*\mathcal{I})^\perp \subseteq j^*\mathcal{J}$.
Let $\beta'':A''\rightarrow B''$ be a morphism in $(j^*\mathcal{I})^\perp$.
Note that we have an isomorphism $\beta'' \cong j^*j_*(\beta'')$.
Thus, to show $\beta'' \in j^*\mathcal{J}$, it suffices to show $j_*(\beta'') \in \mathcal{J}$.
Since $\mathcal{J} = \mathcal{I}^\perp$, we verify this by testing against any $\alpha \in \mathcal{I}$ and any phantom $\phi\in \Phi(\mathfrak{C})$.
Consider the composition $\xi = j_*(\beta'')[1] \circ \phi \circ \alpha$.
Using the adjoint pair ($j^*, j_*$), the morphism $\xi$ corresponds to $\xi'' = \beta''[1] \circ \epsilon_{A''}  \circ j^*(\phi) \circ j^*(\alpha)$, where $\epsilon_{A''}: j^*j_{*}A'' \to A''$ is the counit morphism.
Since $\alpha \in \mathcal{I}$, we have $j^*(\alpha) \in j^*\mathcal{I}$.
Since $\beta'' \in (j^*\mathcal{I})^\perp$, we have $\xi'' = 0$.
Therefore, the composition $\xi$ vanishes.
Thus $j_*(\beta'') \in \mathcal{I}^\perp = \mathcal{J}$.

\textbf{(ii)} We show that ${}^\perp(j^*\mathcal{J}) \subseteq j^*\mathcal{I}$.
Let $\alpha'': A'' \to B''$ be a morphism in ${}^\perp(j^*\mathcal{J})$.
Note that we have an isomorphism $\alpha'' \cong j^*j_!(\alpha'')$.
Thus, to show $\alpha'' \in j^*\mathcal{I}$, it suffices to show $j_!(\alpha'') \in \mathcal{I}$.
Since $\mathcal{I} = {}^\perp\mathcal{J}$, we verify this by testing against any $\beta \in \mathcal{J}$ and any phantom $\phi\in \Phi(\mathfrak{C})$.
Consider the composition $\xi  = \beta[1] \circ \phi \circ j_!(\alpha'')$.
Using the adjoint pair $(j_!, j^*)$, the morphism $\xi$ corresponds to $\xi'' = j^*(\beta)[1] \circ j^*(\phi) \circ \eta_{B''} \circ \alpha''$, where $\eta_{B''}: B'' \to j^*j_!B''$ is the unit morphism.
Since $\beta \in \mathcal{J}$, we have $j^*(\beta) \in j^*\mathcal{J}$.
Since $\alpha'' \in {}^\perp(j^*\mathcal{J})$, we have $\xi'' = 0$.
Therefore, the composition $\xi$ vanishes.
Thus $j_!(\alpha'') \in {}^\perp\mathcal{J} = \mathcal{I}$.

Therefore, we can conclude that $(j^*\mathcal{I}, j^*\mathcal{J})$ is an ideal cotorsion pair in $\mathcal{A}''$.

Due to the Lemma \ref{PC}, if $\mathcal{I}$ is precovering in $\mathcal{A}$, then $j^*\mathcal{I}$ is precovering in $\mathcal{A}''$.
Note that every precoving ideal is special.
Therefore, $(j^*\mathcal{I}, j^*\mathcal{J})$ is a complete ideal cotorsion pair in $\mathcal{A}''$.
\endd

\subsection{Gluing from $\mathcal{T}', \mathcal{T}''$to $\mathcal{T}$}\

\begin{lemm}\label{TTF}
	Consider a recollement as in Definition \ref{recollement}.
	$(j_{!}\mathcal{T}^{\rightarrow }, i_{*}\mathcal{T}^{\rightarrow }), (i_{*}\mathcal{T}^{\rightarrow }, j_*\mathcal{T}^{\rightarrow })$ are complete ideal cotorsion pairs in $\mathcal{T}$.
\end{lemm}
\proof
It is a well-known fact that a recollement induces a TTF triple $(\text{Im } j_{!}, \text{Im } i_{*}, \text{Im } j_{*})$. 
This follows from the adjunctions and the canonical triangles associated with the recollement. 
And cotorsion pairs in $\mathcal{T}$ can induce complete ideal cotorsion pairs with resopect to the proper class of all triangles in $\mathcal{T}$.
\endd

\begin{lemm}\label{equal1}
	Consider a recollement as in Definition \ref{recollement}.
	Let $\mathcal{I}'$ be an ideal in $\mathcal{T}'$ and $\mathcal{I}''$ be an ideal in $\mathcal{T}''$. 
	Then 

	$(1)\ \langle i_*\mathcal{I}',j_*\mathcal{T}''^{\rightarrow }\rangle = \{\xi \in \mathcal{T}\mid i^!(\xi)\in \mathcal{I}'\}$;

	$(2)\ \langle i_*\mathcal{T}'^{\rightarrow }, j_*\mathcal{I}''\rangle = \{\xi \in \mathcal{T}\mid j^*(\xi)\in \mathcal{I}''\}$.

	$(3)\ \langle j_!\mathcal{I}'', i_*\mathcal{T}'^{\rightarrow }\rangle = \{\xi \in \mathcal{T}\mid j^*(\xi)\in \mathcal{I}''\}$;

	$(4)\ \langle j_!\mathcal{T}''^{\rightarrow }, i_*\mathcal{I}'\rangle = \{\xi \in \mathcal{T}\mid i^*(\xi)\in \mathcal{I}'\}$.

\end{lemm}
\proof
We only prove (1); the other cases are similar.

\textbf{Step 1: $\langle i_*\mathcal{I}',j_*\mathcal{T}''^{\rightarrow }\rangle \subseteq \{\xi\in \mathcal{T} \mid i^!(\xi)\in \mathcal{I}'\}$}

Let $\xi: F \to K$ be a morphism in the Toda bracket $\langle i_*\mathcal{I}', j_*\mathcal{T}''^{\rightarrow } \rangle$.
By definition, there exists a triangle $ A \xrightarrow{a} C \xrightarrow{b} B \rightarrow  A[1]$ and morphisms $\xi': F \to C$, $\xi'': C \to K$ such that $\xi = \xi'' \circ \xi'$, $f := \xi'' \circ a \in i_*\mathcal{I}'$ and $g := b \circ \xi' \in j_*\mathcal{T}''$.
\begin{center}
	\begin{tikzcd}
		& & & & F \arrow[d, "g"] \arrow[lld, "\xi'"']\\
		A \arrow[rr, "a"] \arrow[d, "f"'] & & C \arrow[rr, "b"] \arrow[lld, "\xi''"] & & B \\
		K& & & &
	\end{tikzcd}
\end{center}
We need to show $i^!(\xi) \in \mathcal{I}'$.
Apply the functor $i^!$ to the triangle, we get an triangle in $\mathcal{T}'$:
$$ i^!A \xrightarrow{i^!(a)} i^!C \xrightarrow{i^!(b)} i^!B\rightarrow  i^!A[1]. $$
Since $g\in j_{*}\mathcal{T}''$, $g$ factors through an object in the image of $j_*$.
Since $i^!j_* = 0$, we have $i^!(g) = 0$.
Thus, $i^!(b) \circ i^!(\xi') = i^!(g) = 0$.
There exists a morphism $h': i^!F \to i^!A$ such that $i^!(\xi') = i^!(a) \circ h'$.
Substitute this into $i^!(\xi)$: $ i^!(\xi) = i^!(\xi'') \circ i^!(\xi') = i^!(\xi'') \circ i^!(a) \circ h' = i^!(f)\circ h'$.
Note that $f \in i_*\mathcal{I}'$, which implies $i^!(f) \in \mathcal{I}'$ (via the isomorphism $i^!i_* \cong \mathrm{Id}$).
Since $\mathcal{I}'$ is an ideal, the composition $i^!(f)\circ h'$ remains in $\mathcal{I}'$.
Thus, $i^!(\xi) \in \mathcal{I}'$.

\textbf{Step 2: $\{\xi\in \mathcal{T} \mid i^!(\xi)\in \mathcal{I}'\} \subseteq \langle i_*\mathcal{I}',j_*\mathcal{T}''^{\rightarrow}\rangle$}

Let $\xi: F \to K$ be a morphism such that $i^!(\xi) \in \mathcal{I}'$.
Consider the canonical gluing triangle for the object $F$:
$$ i_*i^!F \xrightarrow{\epsilon_F} F \xrightarrow{\eta_F} j_*j^*F \xrightarrow{\delta} i_*i^!F[1] $$
(Here we use $\epsilon_F$ and $\eta_F$ as the counit and unit maps).
We construct a Toda bracket diagram for $\xi$.
\begin{center}
	\begin{tikzcd}
		& && & F \arrow[dd, "\eta_F"] \arrow[lldd,equal] \\
		& && &\\
		i_*i^!F \arrow[rr, "\epsilon_F"] \arrow[d, "i_*i^!(\xi)"'] & & F \arrow[rr, "\eta_F"] \arrow[lldd, "\xi"] & & j_*j^*F\\
		i_*i^!K \arrow[d, "\epsilon_K"'] & && &\\
		K & && &
	\end{tikzcd}
\end{center}
Now we check the conditions for the class memberships.
The compositions $\xi \circ \epsilon_{F}$ and $\epsilon_{K}\circ i_{*}i^!(\xi)$ are both coresponding to $i^!(\xi)$ via the adjoint pair $(i_{*}, i^!)$, that is $\xi \circ \epsilon_{F} = \epsilon_{K}\circ i_{*}i^!(\xi)$.
Since $F\in \mathcal{T}$ and $j^*\mathcal{T}\subseteq \mathcal{T}''$, we have $j^*F\in \mathcal{T}''$.
Then $\eta_{F}=j_{*}(1_{j^*F})\circ \eta_{F}\in j_{*}\mathcal{T}''\subseteq \mathcal{T}$.
Due to $i^!(\xi) \in \mathcal{I}'$, we have $i_{*}i^!(\xi) \in i_*\mathcal{I}'$.
Since $i_{*}\mathcal{I}'$ is an ideal in $\mathcal{T}$, the composition $\epsilon_{K}\circ i_{*}i^!(\xi)$ belongs to $i_*\mathcal{I}'$.
Thus, $\xi \in \langle i_*\mathcal{I}', j_*\mathcal{T}''^{\rightarrow } \rangle$.
\endd

\begin{lemm}\label{perp1}
	Consider a recollement as in Definition \ref{recollement}.
	Let $\mathcal{I}'$ be an ideal in $\mathcal{T}'$ and $\mathcal{I}''$ be an ideal in $\mathcal{T}''$.
	Then 
	
	$(1)\ (i_*\mathcal{I}')^\perp = \{\beta  \in \mathcal{T} \mid i^!(\beta) \in \mathcal{I}'^{\perp}\}$;

	$(2)\ (j_!\mathcal{I}'')^\perp = \{\beta \in \mathcal{T} \mid j^*(\beta) \in (\mathcal{I}'')^{\perp}\}$;

	$(3)\ ^\perp(i_*\mathcal{I}') = \{\alpha \in \mathcal{T} \mid i^*(\alpha) \in {}^{\perp}\mathcal{I}'\}$;
	
	$(4)\ ^\perp(j_*\mathcal{I}'') = \{\alpha \in \mathcal{T} \mid j^*(\alpha) \in {}^{\perp}\mathcal{I}''\}$.
\end{lemm}
\proof
We only prove (1); the other cases are similar.

\textbf{Step 1: $\{\beta \mid i^!(\beta)\in \mathcal{I}'^{\perp}\} \subseteq (i_*\mathcal{I}')^\perp$}

Let $\beta  \in \{\beta \mid i^!(\beta) \in \mathcal{I}'^{\perp}\}$.
We need to show that for any $\alpha  \in i_*\mathcal{I}'$ and any $\phi \in \Phi(\mathfrak{C})$, the composition $\beta [1] \circ \phi \circ \alpha = 0$.
By definition, $\alpha$ factors as $\alpha = f\circ i_*(\alpha')\circ g$ for some $\alpha': A'\rightarrow B' \in \mathcal{I}'$ and morphisms $f, g$ in $\mathcal{T}$.
Consider the composition $\xi = \beta [1] \circ \phi \circ f \circ i_*(\alpha')$.
Using the adjoint pair $(i_{*}, i^!)$, the morphism $\xi$ corresponds to the composition $\xi' = i^!(\beta) [1] \circ i^!(\phi) \circ i^!(f) \circ \eta_{B'} \circ \alpha'$, where $\eta$ is the unit of the adjunction.
Since $i^!(\beta) \in \mathcal{I}'^{\perp}$, $\alpha' \in \mathcal{I}'$ and $i^!(\phi) \in \Phi'$, the composition $\xi'$ vanishes.
Thus, $\xi = 0$.
Hence $\beta  \in (i_*\mathcal{I}')^{\perp}$.

\textbf{Step 2: $(i_*\mathcal{I}')^{\perp} \subseteq \{\beta  \mid i^!(\beta)\in (\mathcal{I}')^{\perp}\}$}

Let $\beta: A\rightarrow B \in (i_*\mathcal{I}')^{\perp}$.
We need to show that $i^!(\beta) \in (\mathcal{I}')^{\perp}$.
Let $\alpha '$ be any morphism in $\mathcal{I}'$ and $\phi'$ be any phantom morphism in $\Phi'$.
It suffices to show that the composition $\xi' = i^!(\beta)[1] \circ \phi' \circ \alpha '$ vanishes.
Using the adjoint pair $(i_*, i^!)$, the morphism $\xi'$ corresponds to the composition $\xi = \beta [1]\circ \epsilon_{A[1]} \circ i_*(\phi') \circ i_*(\alpha ')$ in $\mathcal{T}$, where $\epsilon$ is the counit of the adjunction.
Since $i_*(\alpha ') \in i_*\mathcal{I}'$ and $i_*(\phi') \in \Phi$, the composition $\xi $ vanishes.
Thus, $\xi' = 0$ for any $\alpha' \in \mathcal{I}'$ and $\phi' \in \Phi'$.
That is, $i^!(\beta) \in (\mathcal{I}')^{\perp}$.

Therefore, $(i_*\mathcal{I}')^{\perp} = \{\beta \in \mathcal{T} \mid i^!(\beta) \in (\mathcal{I}')^{\perp}\}$.
\endd

\begin{prop}\label{pair1}
	Consider a recollement as in Definition \ref{recollement}.
	If $(\mathcal{I}',\mathcal{J}')$ is a complete ideal cotorsion pair in $\mathcal{T}'$, then $(i_*\mathcal{I}', \langle i_*\mathcal{J}',j_*\mathcal{T}''\rangle )$ and $(\langle j_!\mathcal{T}'', i_*\mathcal{I}'\rangle ,i_*\mathcal{J}')$  are complete ideal cotorsion pairs in $\mathcal{T}$.
	If $(\mathcal{I}'',\mathcal{J}'')$ is a complete ideal cotorsion pair in $\mathcal{T}''$, then $(j_!\mathcal{I}'', \langle i_*\mathcal{T}',j_*\mathcal{J}''\rangle )$ and $(\langle j_!\mathcal{I}'', i_*\mathcal{T}'\rangle ,j_*\mathcal{J}'')$ are complete ideal cotorsion pairs in $\mathcal{T}$.

\end{prop}
\proof
We only show that $(i_*\mathcal{I}', \langle i_*\mathcal{J}', j_*\mathcal{T}''\rangle)$ is a complete ideal cotorsion pair; the other cases are similar.
$\mathcal{I}'$ is a special precovering ideal in $\mathcal{T}'$.
Then $i_*\mathcal{I}'$ is a special precovering ideal in $\mathcal{T}$.
Thus $(i_{*}\mathcal{I}', (i_{*}\mathcal{I}')^\perp)$ is a complete ideal cotorsion pair in $\mathcal{T}$.
by Lemma\ref{equal1} and Lemma\ref{perp1}, we have $(i_*\mathcal{I}')^\perp = \langle i_*\mathcal{J}', j_*\mathcal{T}''\rangle$.
Therefore, $(i_*\mathcal{I}', \langle i_*\mathcal{J}', j_*\mathcal{T}''\rangle)$ is a complete ideal cotorsion pair in $\mathcal{T}$.
\endd

\begin{lemm}
	Consider a recollement as in Definition \ref{recollement}.
	Let $\mathcal{I}'$ be an ideal in $\mathcal{T}'$ and $\mathcal{I}''$ be an ideal in $\mathcal{T}''$.
	Define a subset $\mathcal{I}$ of $\mathcal{T}$ by
	$$ \mathcal{I} = \{\alpha \in \mathcal{T}\mid i^*(\alpha) \in \mathcal{I}', \ j^*(\alpha) \in \mathcal{I}'' \}. $$
	Then $\mathcal{I}$ is an ideal in $\mathcal{T}$. Moreover, the following equalities hold:

	$(1)\ i^*\mathcal{I} = \mathcal{I}' = \{\alpha  \in \mathcal{T}' \mid i_*(\alpha) \in \mathcal{I}\}$;

	$(2)\ j^*\mathcal{I} = \mathcal{I}'' = \{\alpha  \in \mathcal{T}''\mid j_!(\alpha) \in \mathcal{I}\}$.
\end{lemm}

\proof
Since $i^*$ and $j^*$ are additive functors and $\mathcal{I}', \mathcal{I}''$ are ideals, it is clear that $\mathcal{I}$ is closed under addition and composition. Thus, $\mathcal{I}$ is an ideal in $\mathcal{T}$.

By definition, for any $\alpha  \in \mathcal{I}$, $i^*(\alpha) \in \mathcal{I}'$, so $i^*\mathcal{I} \subseteq \mathcal{I}'$.
Conversely, let $\alpha  \in \mathcal{I}'$.
Since $j^*i_* = 0$ implies $j^*(i_*(\alpha )) = 0 \in \mathcal{I}''$, and $i^*(i_*(\alpha)) \cong \alpha \in \mathcal{I}'$, we have $i_*(\alpha) \in \mathcal{I}$.
That is $\mathcal{I}'\subseteq \{\alpha  \in \mathcal{T}'\mid i_{*}(\alpha)\in \mathcal{I}\}$.
Thus, $\alpha \cong i^*(i_*(\alpha)) \in i^*\mathcal{I}$. This shows $i^*\mathcal{I} = \mathcal{I}'$.
If $\alpha  \in \{\alpha \in \mathcal{T}'\mid i_{*}(\alpha)\in \mathcal{I}\}$, then $i_*(\alpha) \in \mathcal{I}$.
By the definition of $\mathcal{I}$, this implies $i^*(i_*(\alpha)) \in \mathcal{I}'$.
Since $i^*i_* \cong \text{Id}_{\mathcal{T}'}$, we have $\alpha \in \mathcal{I}'$, so $\{\alpha  \in \mathcal{T}'\mid i_*(\alpha)\in \mathcal{I}\} \subseteq \mathcal{I}'$.

Similarly, by definition, $j^*\mathcal{I} \subseteq \mathcal{I}''$.
Conversely, let $\alpha \in \mathcal{I}''$.
Since $i^*j_! = 0$ implies $i^*(j_!(\alpha)) = 0 \in \mathcal{I}'$, and $j^*(j_!(\alpha)) \cong \alpha \in \mathcal{I}''$, we have $j_!(\alpha) \in \mathcal{I}$.
This means $\mathcal{I}'' \subseteq \{\alpha \in \mathcal{T}''\mid j_!(\alpha) \in \mathcal{I}\}$.
Thus, $\alpha \cong j^*(j_!(\alpha)) \in j^*\mathcal{I}$. This shows $j^*\mathcal{I} = \mathcal{I}''$.
If $\alpha \in \{\alpha \in \mathcal{T}''\mid j_!(\alpha) \in \mathcal{I}\}$, then $j_!(\alpha) \in \mathcal{I}$. By the definition of $\mathcal{I}$, $j^*(j_!(\alpha)) \in \mathcal{I}''$.
Since $j^*j_! \cong \text{Id}_{\mathcal{T}''}$, we have $\alpha \in \mathcal{I}''$, so $\{\alpha \in \mathcal{T}''\mid j_!(\alpha) \in \mathcal{I}\} \subseteq \mathcal{I}''$.
\endd

In the same way, we claim following conclusion without proof.
\begin{lemm}
	Let $\mathcal{J}'$ be an ideal in $\mathcal{T}'$, $\mathcal{J}''$ be an ideal in $\mathcal{T}''$.
	Define a subset $\mathcal{J}$ of $\mathcal{T}$ by
	$$ \mathcal{J} = \{\beta  \in \mathcal{T}\mid i^!(\beta) \in \mathcal{J}', \ j^*(\beta) \in \mathcal{J}'' \}. $$
	Then $\mathcal{J}$ is an ideal in $\mathcal{T}$. Moreover, the following equalities hold:

	$(1)\ i^!\mathcal{J} = \mathcal{J}' = \{\beta  \in \mathcal{T}'\mid i_*(\beta) \in \mathcal{J}\}$;

	$(2)\ j^*\mathcal{J} = \mathcal{J}'' = \{\beta  \in \mathcal{T}''\mid j_*(\beta) \in \mathcal{J}\}$.
\end{lemm}

\noindent \textit{\textbf{Proof of Theorem \ref{main1}}}:

Assume that $(\mathcal{I}',\mathcal{J}')$ is a complete ideal cotorsion pair in $\mathcal{T}'$ and $(\mathcal{I}'',\mathcal{J}'')$ is a complete ideal cotorsion pair in $\mathcal{T}''$.
According to Proposition \ref{pair1}, we have complete ideal cotorsion pairs $(i_*\mathcal{I}', \langle i_*\mathcal{J}',j_*\mathcal{T}''\rangle)$ and $(j_!\mathcal{I}'', \langle i_*\mathcal{T}',j_*\mathcal{J}''\rangle)$ in $\mathcal{T}$.
In addition, Lemma \ref{Toda} guarantees the existence of a complete ideal cotorsion pair
\[
	(\langle j_!\mathcal{I}'',i_*\mathcal{I}'\rangle, \langle i_*\mathcal{J}',j_*\mathcal{T}''\rangle \langle i_*\mathcal{T}',j_*\mathcal{J}''\rangle).
\]

Next, we claim that the product of ideals $\langle i_*\mathcal{J}',j_*\mathcal{T}''\rangle \langle i_*\mathcal{T}',j_*\mathcal{J}''\rangle$ is equal to $\langle i_*\mathcal{J}',j_*\mathcal{J}''\rangle$.

First, let $\phi$ be any morphism in $\langle i_*\mathcal{J}',j_*\mathcal{T}''\rangle \langle i_*\mathcal{T}',j_*\mathcal{J}''\rangle$.
By definition, there exists a commutative diagram as follows:
\begin{center}
	\begin{tikzcd}
		& & & && & & & H \arrow[d, "\beta_2"] \arrow[lld,"\xi_1"'] \\
		& & & & D \arrow[d, "f_1"] \arrow[rr,"d"] & & E \arrow[rr,"f"'] \arrow[lld,"\xi_2"] & & F \\
		& & & & K \arrow[d, "f_2"] \arrow[lld,"\xi_3"'] & & & &\\
		C \arrow[rr] \arrow[d, "\beta_1"] & & B \arrow[rr] \arrow[lld,"\xi_4"] & &  A & & & &\\
		M & & & && & & &
	\end{tikzcd}
\end{center}
where $\beta_1 \in i_*\mathcal{J}'$, $f_2 \in j_*\mathcal{T}''$, $f_1 \in i_*\mathcal{T}'$, $\beta_2 \in j_*\mathcal{J}''$.
The morphism $\phi$ is given by the composition $\phi = \xi_4 \circ \xi_3 \circ \xi_2 \circ \xi_1$.
Consider the composition $f_2 \circ f_1$. Since $f_1$ factors through $i_*(\mathcal{T}')$ and $f_2$ factors through $j_*(\mathcal{T}'')$, and knowing that $j^*i_* = 0$, the orthogonality of the recollement implies that:
\[
	\text{Hom}_{\mathcal{T}}(i_*X', j_*Y'') \cong \text{Hom}_{\mathcal{T}''}(j^*i_*X', Y'') = 0.
\]
Thus, we have $f_2 \circ f_1 = 0$.
Consequently, $f_2 \circ \xi_2 \circ d = f_2 \circ f_1 = 0$.
The vanishing of this composition implies the existence of a morphism $\delta: F \to A$ such that $f_2 \circ \xi_2 = \delta \circ f$.
This yields the following commutative diagram:
\begin{center}
	\begin{tikzcd}
		& & & && & & & H \arrow[d, "\beta_2"] \arrow[lld,"\xi_1"'] \\
		& & & & D \arrow[d, "f_1"] \arrow[rr,"d"] & & E \arrow[rr,"f"'] \arrow[lld,"\xi_2"] & & F \arrow[lllldd, "\delta", dashed] \\
		& & & & K \arrow[d, "f_2"] \arrow[lld,"\xi_3"'] & & & &\\
		C \arrow[rr] \arrow[d, "\beta_1"] & & B \arrow[rr] \arrow[lld,"\xi_4"] & &  A & & & &\\
		M & & & && & & &
	\end{tikzcd}
\end{center}
Since $\delta \circ \beta_2 \in j_*\mathcal{J}''$ (as $\mathcal{J}''$ is an ideal), we conclude that $\phi \in \langle i_*\mathcal{J}',j_*\mathcal{J}''\rangle$.

On the other hand, let $\phi$ be any morphism in $\langle i_*\mathcal{J}',j_*\mathcal{J}''\rangle$.
Then there exists a commutative diagram:

\begin{center}
	\begin{tikzcd}
		& & & & N \arrow[d, "d\circ j_*(\beta'')\circ c"] \arrow[lld, "\xi_1"'] \\
		E \arrow[d, "b\circ i_*(\beta')\circ a"'] \arrow[rr, "e"] & & F \arrow[rr, "f"] \arrow[lld, "\xi_2"] & & G\\
		M& & & &
	\end{tikzcd}
\end{center}
with $\beta' \in \mathcal{J}'$, $\beta'' \in \mathcal{J}''$, and $\phi = \xi_2 \circ \xi_1$.
Taking the cobase change along $a$ and the base change along $d$, we obtain the following commutative diagram:
\begin{center}
	\begin{tikzcd}
		& & N \arrow[rd, "c"] \arrow[rddd, "\xi_1" description] \arrow[dd, "n" description, dotted] & && \\
		& & & j_*C'' \arrow[rd, "j_*(\beta'')"] && \\
		i_*A' \arrow[dd,equal] \arrow[rr,"a_1"description] & & L \arrow[dd, "l" description,near start] \arrow[rr,"d_2"description,near end]& & j_*D'' \arrow[dd, "d" description,near start] & \\
		& E \arrow[rr, "e" description,near start] \arrow[ld, "a" description] & & F \arrow[rr, "f" description,near start] \arrow[ld, "h" description] \arrow[lddd, "\xi_2" description] && G \arrow[ld,equal] \\
		i_*A' \arrow[rr] \arrow[rd, "i_*(\beta')" '] & & H \arrow[rr,"g"description,near end] \arrow[dd, "m" description, dotted]& & G & \\
		& i_*B' \arrow[rd, "b"'] & & && \\
		& & M & &&
	\end{tikzcd}
\end{center}
Since the diagram commutes, specifically $\xi_2 \circ e = b \circ i_*(\beta') \circ a$, there exists a morphism $m: H \to M$ such that $\xi_2 = m \circ h$.
Similarly, $d \circ j_*(\beta'') \circ c = f \circ \xi_1 = g \circ h \circ \xi_1$ implies that there exists a morphism $n: N \to L$ such that $h \circ \xi_1 = l \circ n$.
Thus, we have:
\[
	\phi = \xi_2 \circ \xi_1 = m \circ l \circ n.
\]
This factorization implies that $\phi \in \langle i_*\mathcal{J}',j_*\mathcal{T}''\rangle \langle i_*\mathcal{T}',j_*\mathcal{J}''\rangle$.
Therefore, we have established the equality:
\[
	\langle i_*\mathcal{J}',j_*\mathcal{T}''\rangle \langle i_*\mathcal{T}',j_*\mathcal{J}''\rangle = \langle i_*\mathcal{J}',j_*\mathcal{J}''\rangle.
\]
It follows that $(\langle j_!\mathcal{I}'',i_*\mathcal{I}'\rangle, \langle i_*\mathcal{J}',j_*\mathcal{J}''\rangle)$ is a complete ideal cotorsion pair in $\mathcal{T}$.

\medskip

Conversely, assume that $(\langle j_!\mathcal{I}'',i_*\mathcal{I}'\rangle, \langle i_*\mathcal{J}',j_*\mathcal{J}''\rangle)$ is a complete ideal cotorsion pair in $\mathcal{T}$.
It follows from Lemma \ref{main2} that the pairs obtained by applying the functors $i^*, i^!$ and $j^*$ preserve the structure.
Specifically:
\begin{enumerate}
	\item $(i^*\langle j_!\mathcal{I}'',i_*\mathcal{I}'\rangle, i^!\langle i_*\mathcal{J}',j_*\mathcal{J}''\rangle)$ is a pair in $\mathcal{T}'$.
	\item $(j^*\langle j_!\mathcal{I}'',i_*\mathcal{I}'\rangle, j^*\langle i_*\mathcal{J}',j_*\mathcal{J}''\rangle)$ is a pair in $\mathcal{T}''$.
\end{enumerate}
Using the properties of Toda brackets and the functors (noting that $i^*i_* \cong \text{Id}$, $j^*j_! \cong \text{Id}$, etc., and the orthogonality relations), we calculate:
\begin{align*}
	i^*\langle j_!\mathcal{I}'',i_*\mathcal{I}'\rangle & = \langle i^*j_!\mathcal{I}'',i^*i_*\mathcal{I}'\rangle = \langle 0, \mathcal{I}'\rangle = \mathcal{I}';   \\
	i^!\langle i_*\mathcal{J}',j_*\mathcal{J}''\rangle & = \langle i^!i_*\mathcal{J}',i^!j_*\mathcal{J}''\rangle = \langle \mathcal{J}', 0\rangle = \mathcal{J}';   \\
	j^*\langle j_!\mathcal{I}'',i_*\mathcal{I}'\rangle & = \langle j^*j_!\mathcal{I}'',j^*i_*\mathcal{I}'\rangle = \langle \mathcal{I}'', 0\rangle = \mathcal{I}''; \\
	j^*\langle i_*\mathcal{J}',j_*\mathcal{J}''\rangle & = \langle j^*i_*\mathcal{J}',j^*j_*\mathcal{J}''\rangle = \langle 0, \mathcal{J}''\rangle = \mathcal{J}''.
\end{align*}
Therefore, $(\mathcal{I}',\mathcal{J}')$ and $(\mathcal{I}'',\mathcal{J}'')$ are complete ideal cotorsion pairs.
\qed

\noindent \emph{\textbf{Remark}}:
By using Toda bracket and intersection, there are many ideal cotorsion pairs.
In order to show the relationship between different ideal cotorsion pairs constructed from $(\mathcal{I}',\mathcal{J}')$ and $(\mathcal{I}'',\mathcal{J}'')$ more directly, we have the diagram in the appendix.

\newpage
\bibliography{ideal}

\section{Appendix}

\begin{center}
	\begin{figure}[!ht]
		\centering
		\includegraphics[width=0.8\linewidth]{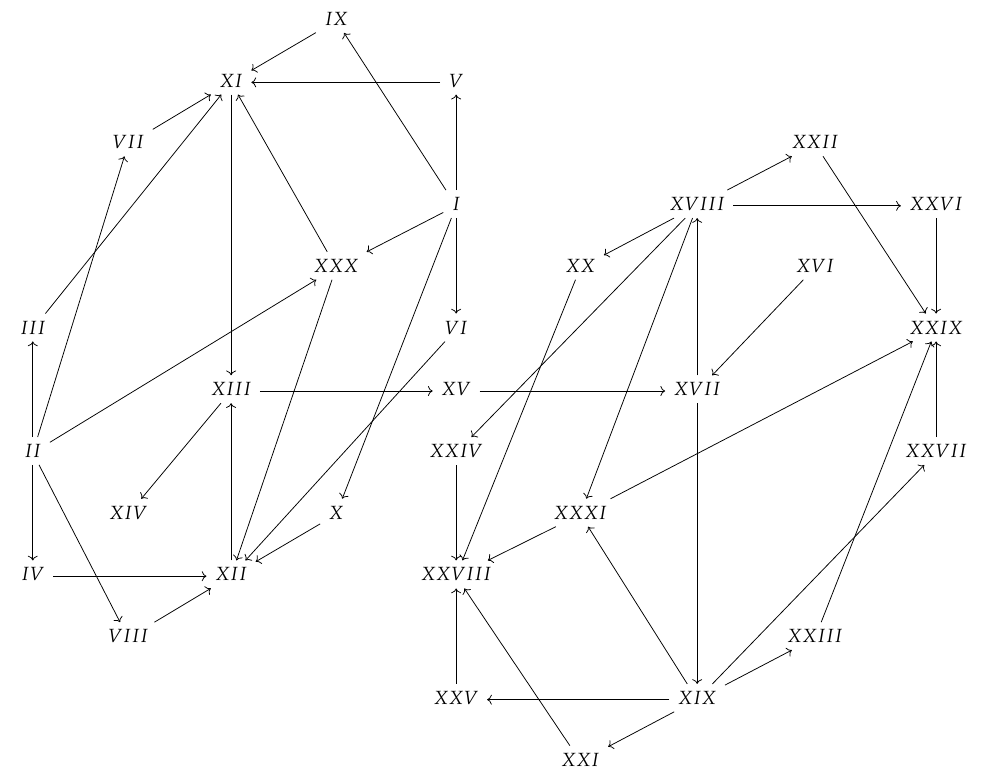}
	\end{figure}
\end{center}

\begin{center}
	\begin{figure}[!ht]
		\centering
		\includegraphics[width=\linewidth]{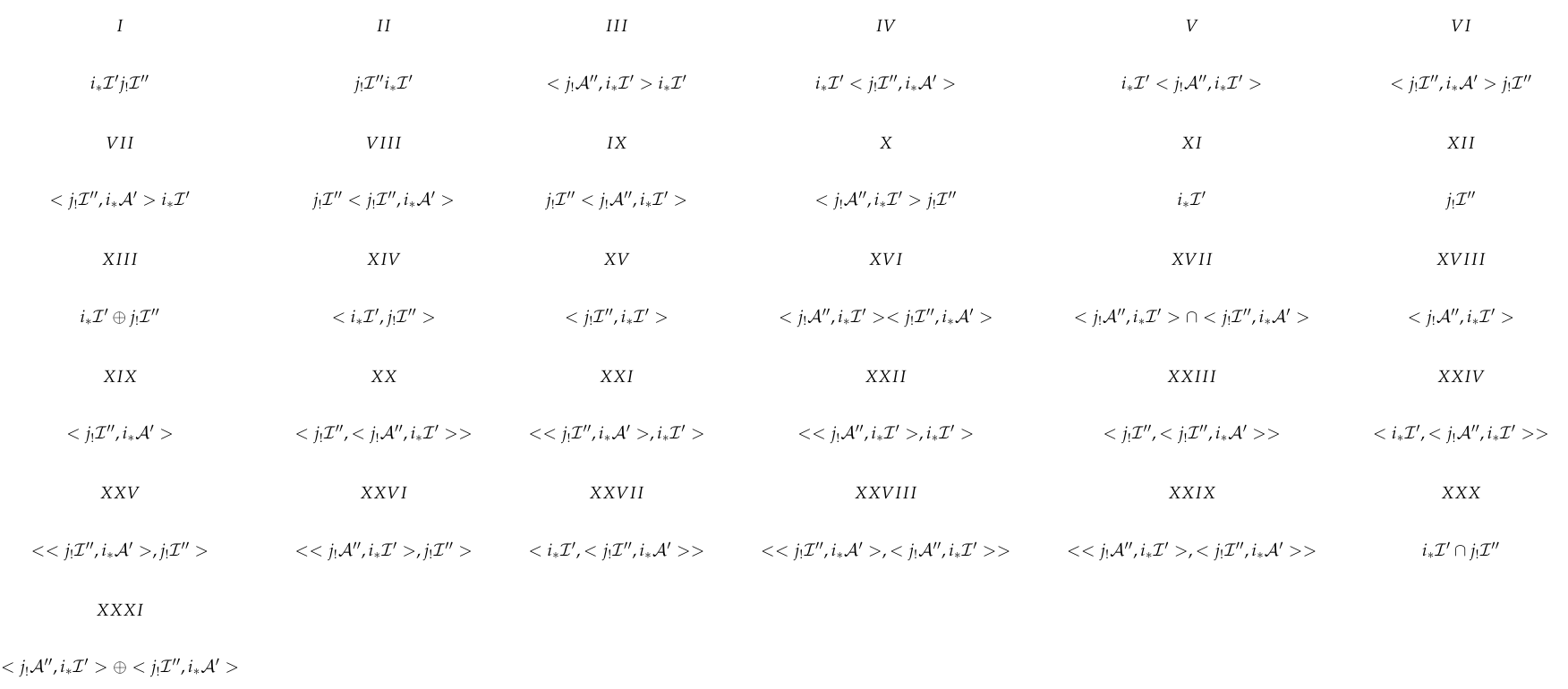}
	\end{figure}
\end{center}
The arrows in the above diagram represent embedded maps. All ideals in the diagram are ideal cotorsion pair with their right orthogonal class.

\vspace{4mm}

%Authors
%\noindent\textbf{Hongxing Chen}\\
%School of Mathematical Sciences \&  Academy for Multidisciplinary Studies, \\Capital Normal University, Beijing 100048, P. R. China;\\
%Email: \textsf{chenhx@cnu.edu.cn}\\[1mm]
%\textbf{Jiangsheng Hu}\\
%School of Mathematics, Hangzhou Normal University, Hangzhou 311121, P. R. China.\\
%Email: \textsf{jiangshenghu@hotmail.com}\\[1mm]
\end{document}